\documentclass{svmult}
\usepackage{amscd,amssymb,amsmath}

\usepackage[colorlinks,plainpages,urlcolor=blue]{hyperref}

\newcommand{\MR}[1]
{\href{http://www.ams.org/mathscinet-getitem?mr=#1}{MR#1}}
\newcommand{\MRh}[2]
{\href{http://www.ams.org/mathscinet-getitem?mr=#1}{MR#1 (#2)}}

\newcommand{\C}{{\mathbb C}}  
\newcommand{\Z}{{\mathbb Z}} 
\newcommand{\N}{{\mathbb N}} 
\renewcommand{\k}{\Bbbk} 

\DeclareMathOperator{\Hom}{Hom}
\DeclareMathOperator{\im}{im}

\newcommand{\DSA}{{\mathcal A}} 
\newcommand{\DSB}{{\mathcal B}}
\newcommand{\DSV}{{\mathcal V}} 
\newcommand{\DSX}{{\mathcal X}} 

\begin{document}

\title*{Torsion in the homology of {M}ilnor fibers of hyperplane arrangements}
\titlerunning{{M}ilnor fibrations of arrangements}

\author{Graham Denham and Alexander~I.~Suciu}
\institute{Graham Denham \at 
Department of Mathematics, University of Western Ontario,
London, ON  N6A 5B7\\ 
\email{gdenham@uwo.ca}\\
Supported by NSERC
\and
Alexander~I.~Suciu \at 
Department of Mathematics,
Northeastern University,
Boston, MA 02115\\ 
\email{a.suciu@neu.edu}\\
Supported in part by NSF grant DMS--1010298 and NSA grant H98230-13-1-0225
}

\maketitle

\begin{abstract}
{As is well-known, the homology groups of the complement 
of a complex hyperplane arrangement are torsion-free.  
Nevertheless, as we showed in a recent paper \cite{DS_DS13}, 
the homology groups of the Milnor fiber of such an arrangement 
can have non-trivial integer torsion. We give here a brief account of the 
techniques that go into proving this result, outline some of 
its applications, and indicate some further questions that
it brings to light.}
\end{abstract}

\subsubsection*{Introduction}
\label{intro}

This talk reported on the main results of \cite{DS_DS13}.  Here, we give
an outline of our approach and a summary of our conclusions.  Our
main result gives a construction of a family of projective hypersurfaces
for which the Milnor fiber has torsion in homology.  The hypersurfaces
we use are hyperplane arrangements, for which techniques are available
to examine the homology of finite cyclic covers quite explicitly, by
reducing to rank $1$ local systems.

The parameter spaces for rank $1$ local systems with non-vanishing 
homology are known as characteristic varieties.  In the special case of 
complex hyperplane arrangement complements, the 
combinatorial theory of multinets largely elucidates their structure, at 
least in degree $1$.  We make use of an iterated parallel connection 
construction to build arrangements with suitable characteristic varieties, 
then vary the characteristic of the field of definition in order to construct 
finite cyclic covers with torsion in first homology.  These covers include 
the Milnor fiber.  We now give some detail about each step.

\subsubsection*{The Milnor fibration}
\label{mf}

A classical construction due to J.~Milnor associates 
to every homogeneous polynomial $f\in \C[z_1,\dots , z_\ell]$ 
a fiber bundle, with base space $\C^*=\C\setminus \{0\}$, 
total space the complement in $\C^{\ell}$ to the hypersurface 
defined by $f$, and projection map 
$f\colon \C^{\ell} \setminus f^{-1}(0)\to \C^*$.

The Milnor fiber $F=f^{-1}(1)$ has the homotopy type of a finite,
$(\ell-1)$-dimensional CW-complex, while the monodromy of the 
fibration, $h\colon F\to F$, is given by $h(z)=e^{2\pi i/n} z$, 
where $n$ is the degree of $f$. 
If $f$ has an isolated singularity at the origin, then $F$ is 
homotopic to a bouquet of $(\ell-1)$-spheres, 
whose number can be determined by algebraic means. 
In general, though, it is a rather hard problem to compute 
the homology groups of the Milnor fiber, even in the case when 
$f$ completely factors into distinct linear forms: that is, 
when the hypersurface $\{f=0\}$ is a hyperplane arrangement.

Building on our previous work with D.~Cohen \cite{DS_CDS03}, 
we show there exist projective hypersurfaces (indeed, hyperplane 
arrangements) whose complements have torsion-free 
homology, but whose Milnor fibers have torsion in homology. 
Our main result can be summarized as follows. 

\begin{theorem}
\label{thm:main}
For every prime $p\ge 2$, there is a hyperplane arrangement 
whose Milnor fiber has non-trivial $p$-torsion in homology. 
\end{theorem}
  
This resolves a problem posed by Randell~\cite[Problem 7]{DS_Ra11},
who conjectured that Milnor fibers of hyperplane arrangements 
have torsion-free homology.  Our examples also give a 
refined answer to a question posed by Dimca and 
N\'emethi~\cite[Question~3.10]{DS_DN04}: torsion in homology 
may appear even if the hypersurface is defined by a reduced 
equation.  We note the following consequence:

\begin{corollary}
There are hyperplane arrangements whose Milnor fibers do not
have a minimal cell structure.
\end{corollary}

This stands in contrast with arrangement complements, which
always admit perfect Morse functions. 
Our method also allows us to compute the homomorphism induced 
in homology by the monodromy, with coefficients 
in a field of characteristic not dividing the order of the monodromy.  

It should be noted that our approach produces only 
examples of arrangements $\DSA$ for which the Milnor fiber $F(\DSA)$ 
has torsion in $q$-th homology, for some $q>1$.  This leaves 
open the following question.

\begin{question}
\label{q:h1}
Is there an arrangement $\DSA$ such that $H_1(F(\DSA),\Z)$ has 
non-zero torsion? 
\end{question}

Since our methods rely on complete reducibility, it is also natural to
ask: do there exist projective
hypersurfaces of degree $n$ for which the Milnor fiber has homology
$p$-torsion, where $p$ divides $n$?  If so, is there a hyperplane arrangement
with this property?

A much-studied question in the subject is whether the Betti 
numbers of the Milnor fiber of an arrangement  $\DSA$ are determined 
by the intersection lattice, $L(\DSA)$. While we do not address this 
question directly, our result raises a related, and arguably 
even more subtle problem.

\begin{question}
\label{q:comb tors}
Is the torsion in the homology of the Milnor fiber of a hyperplane 
arrangement combinatorially determined? 
\end{question}

As a preliminary question, one may also ask: can one predict 
the existence of torsion in the homology of $F(\DSA)$ 
simply by looking at $L(\DSA)$?  As it turns out, under fairly general 
assumptions, the answer is yes: if $L(\DSA)$ satisfies certain very 
precise conditions, then automatically $H_*(F(\DSA),\Z)$ will 
have non-zero torsion, in a combinatorially determined degree.

\subsubsection*{Hyperplane arrangements}
Let $\DSA$ be a (central) arrangement of $n$ hyperplanes in $\C^{\ell}$, 
defined by a polynomial $Q(\DSA)=\prod_{H\in \DSA} f_H$, where each 
$f_H$ is a linear form whose kernel is $H$.   The starting point 
of our study is the well-known observation that the Milnor fiber 
of the arrangement, $F(\DSA)$, is a cyclic, $n$-fold regular cover 
of the projectivized complement, $U(\DSA)$; this cover is 
defined by the homomorphism $\delta\colon \pi_1(U(\DSA))
\twoheadrightarrow \Z_n$, taking each meridian generator $x_H$ to $1$. 

Now, if $\k$ is an algebraically closed field whose characteristic 
does not divide $n$, then $H_q(F(\DSA),\k)$ decomposes as a 
direct sum, $\bigoplus_{\rho} H_q(U(\DSA),\k_{\rho})$, where the 
rank $1$ local systems $\k_{\rho}$ are indexed by characters 
$\rho\colon \pi_1(U(\DSA))\to \k^*$ that factor through $\delta$.  
Thus, if there is such a character $\rho$ for which  
$H_q(U(\DSA),\k_{\rho})\ne 0$, but there is no corresponding 
character in characteristic $0$, then the 
group $H_q(F(\DSA),\Z)$ will have non-trivial $p$-torsion.

To find such characters, we first consider multi-arrangements 
$(\DSA,m)$, with positive integer weights $m_H$ attached 
to each hyperplane $H\in \DSA$.  The corresponding Milnor 
fiber, $F(\DSA,m)$, is defined by the homomorphism 
$\delta_m\colon \pi_1(U(\DSA))\twoheadrightarrow \Z_N$, $x_H\mapsto m_H$, 
where $N$ denotes the sum of the weights. Fix a prime $p$. 
Starting with an arrangement $\DSA$ supporting a suitable multinet, 
we find a deletion $\DSA' =\DSA\setminus \{H\}$, and  
a choice of multiplicities $m'$ on $\DSA'$ such that $H_1(F(\DSA',m'),\Z)$ 
has $p$-torsion. Finally, we construct a ``polarized" arrangement 
$\DSB=\DSA' \| m'$, and show that $H_*(F(\DSB),\Z)$ has $p$-torsion.

\subsubsection*{Characteristic varieties}
\label{cv}
Our arguments depend on properties of the jump loci of rank $1$ local
systems.  The {\em characteristic varieties}\/ of a connected, finite 
CW-complex $X$ are the subvarieties $\DSV^q_d(X,\k)$ of the character group 
$\widehat{G}=\Hom(G,\k^*)$, consisting of those characters 
$\rho$ for which $H_q(X,\k_{\rho})$ had dimension at least $d$.

Suppose $X^{\chi}\to X$ is a regular cover, defined by an epimorphism 
$\chi\colon G\to A$ to a finite abelian group, and if $\k$ is an algebraically 
closed field of characteristic $p$, where $p\nmid \left|\DSA\right|$, then 
$\dim H_q(X^{\chi},\k)=\sum_{d\geq 1}\left|\im(\widehat{\chi}_\k) \cap  
\DSV_d^q(X,\k)\right|$, where $\widehat{\chi}_\k\colon \widehat{A} \to \widehat{G}$ 
is the induced morphism between character groups.

\begin{theorem}
\label{thm:ptors1}
Let $X^{\chi}\to X$ be a regular, finite cyclic cover. 
Suppose that $\im(\widehat{\chi}_\C)\not\subseteq \DSV^q_1(X,\C)$, but
$\im(\widehat{\chi}_\k) \subseteq\DSV^q_1(X,\k)$, 
for some field $\k$ of characteristic $p$ not dividing the order 
of the cover. Then $H_q(X^{\chi},\Z)$ has non-zero $p$-torsion.  
\end{theorem}

\subsubsection*{Multinets}
\label{multinets}
In the case when $X=M(\DSA)$ is the complement of a hyperplane
arrangement, the positive-dimensional components of
the characteristic variety $\DSV^1_1(X,\C)$ have a combinatorial
description, for which we refer in particular to work of Falk 
and Yuzvinsky in \cite{DS_FY}.  

A {\em multinet}\/ consists of a partition of $\DSA$ into at least $3$ 
subsets $\DSA_1,\ldots,\DSA_k$, together with an assignment of 
multiplicities, $m\colon \DSA\to \N$, and a subset $\DSX$ of the rank 
$2$ flats, such that any two 
hyperplanes from different parts intersect 
at a flat in $\DSX$, and several technical conditions are satisfied:   
for instance, the sum of the multiplicities over each part $\DSA_i$  
is constant, and for each flat $Z\in \DSX$, the sum 
$n_{Z}:=\sum_{H\in\DSA_i\colon H\supset Z} m_H$ is independent of $i$.  
Each multinet gives rise to an orbifold fibration 
$X\to \mathbb{P}^1\setminus \{\text{$k$ points}\}$;  
in turn, such a map yields by pullback 
an irreducible component of $\DSV^1_1(X,\C)$. 

We say that a multinet on $\DSA$ is {\em pointed}\/ if
for some hyperplane $H$, we have $m_{H}>1$ and $m_{H} \mid n_Z$ 
for each flat $Z\subset H$ in $\DSX$.
We show that the complement of the deletion $\DSA':=\DSA\setminus \{H\}$ 
supports an orbifold fibration with base $\C^*$ and at least one 
multiple fiber. Consequently, for any prime $p\mid  m_{H}$, 
and any sufficiently large integer $r$ not divisible by 
$p$, there exists a regular, $r$-fold cyclic cover $Y\to U(\DSA')$ 
such that $H_1(Y,\Z)$ has $p$-torsion.

Furthermore, we also show that any finite cyclic cover 
of an arrangement complement is dominated by a Milnor 
fiber corresponding to a suitable choice of multiplicities. 
Putting things together, we obtain the following result.

\begin{theorem}
\label{thm:multi tors}
Suppose $\DSA$ admits a pointed multinet, with distinguished 
hyperplane $H$ and multiplicity vector $m$.  Let $p$ be a prime 
dividing $m_H$. There is then a choice of multiplicity vector $m'$ 
on the deletion $\DSA' =\DSA\setminus \{H\}$ such that 
$H_1(F(\DSA',m'),\Z)$ has non-zero $p$-torsion.
\end{theorem}

For instance, if $\DSA$ is the reflection arrangement of type 
$\operatorname{B}_3$, defined by the polynomial 
$Q=xyz(x^2-y^2)(x^2-z^2)(y^2-z^2)$, then $\DSA$ 
satisfies the conditions of the theorem, for 
$m=(2,2,2,1,1,1,1,1,1)$ and $H= \{z=0\}$.  
Choosing then multiplicities $m'=(2,1,3,3,2,2,1,1)$ 
on $\DSA'$ shows that $H_1(F(\DSA',m'),\Z)$ has non-zero 
$2$-torsion.  

Similarly, for primes $p>2$, we use the fact that
the reflection arrangement of the full monomial complex reflection
group, $\DSA(p,1,3)$, admits a pointed multinet. 
This yields $p$-torsion in the first homology 
of the Milnor fiber of a suitable multi-arrangement on 
the deletion.   

\subsubsection*{Parallel connections and polarizations}
\label{parallel}

The last step of our construction replaces multi-arrangements with
simple arrangements.  We add more hyperplanes and increase the
rank by means of suitable iterated parallel
connections.  The complement of the parallel connection of 
two arrangements is diffeomorphic to the product of the respective
complements, by a result of Falk and Proudfoot~\cite{DS_FP}. Then the
characteristic varieties of the result are given by a 
formula due to Papadima and Suciu~\cite{DS_PS10}.  

We organize the process by noting that parallel
connection of matroids has an operad structure, and we analyze
a special case which we call the {\em polarization} of a multi-arrangement
$(\DSA,m)$.  By analogy with a construction involving monomial ideals,
we use parallel connection to attach to each hyperplane $H$ the
central arrangement of $m_H$ lines in $\C^2$, to obtain a simple
arrangement we denote by $\DSA\| m$.  
A crucial point here is the connection between the respective 
Milnor fibers:  the pullback of the cover $F(\DSA\| m)\to U(\DSA\| m)$ 
along the canonical inclusion $U(\DSA) \to U(\DSA\| m) $ 
is equivalent to the cover $F(\DSA,m)\to U(\DSA)$.  
Using this fact, we prove the following.

\begin{theorem}
\label{thm:polar tors}
Suppose $\DSA$ admits a pointed multinet, with distinguished 
hyperplane $H$ and multiplicity $m$.  Let $p$ be a prime 
dividing $m_H$. 
There is then a choice of multiplicities $m'$ on the deletion 
$\DSA' =\DSA\setminus \{H\}$ such that the Milnor fiber of the 
polarization $\DSA' \| m'$ has $p$-torsion in homology, 
in degree $1+\left|\{K\in \DSA':  m'_K\ge 3\}\right|$.
\end{theorem}

For instance, if $\DSA'$ is the deleted $\operatorname{B}_3$ 
arrangement as above, then choosing 
$m'=(8,1,3,3,5,5,1,1)$ produces an arrangement 
$\DSB=\DSA'\| m'$ of $27$ hyperplanes in $\C^8$, such that 
$H_6(F(\DSB),\Z)$ has $2$-torsion of rank $108$.


\begin{thebibliography}{00}

\bibitem{DS_CDS03} D.C.~Cohen, G.~Denham, A.I.~Suciu, 
{\em Torsion in {M}ilnor fiber homology}, Alg. Geom. 
Topology \textbf{3} (2003), 511--535. 
\MRh{1997327}{2004d:32043}

\bibitem{DS_DS13} G.~Denham, A.I.~Suciu,
{\em Multinets, parallel connections, and {M}ilnor fibrations of
arrangements}, Proc.~London Math.~Soc. \textbf{108} (2014), 
no.~6, 1435--1470.
\MR{3218315}

\bibitem{DS_DN04} A.~Dimca, A.~N\'emethi,
{\em Hypersurface complements, Alexander modules and monodromy},
in: Real and complex singularities, 19--43,
Contemp. Math., vol.~354, Amer. Math. Soc., Providence, RI, 2004. 
\MRh{2087802}{2005h:32076}

\bibitem{DS_FP} M.~Falk, N.~Proudfoot, 
{\em Parallel connections and bundles of arrangements},  
Topology Appl. \textbf{118} (2002),  no.~1-2, 65--83.   
\MRh{1877716}{2002k:52033}

\bibitem{DS_FY} M.~Falk, S.~Yuzvinsky,
{\em Multinets, resonance varieties, and pencils of plane curves},
Compositio Math. \textbf{143} (2007), no.~4, 1069--1088.
\MRh{2339840}{2009e:52043}

\bibitem{DS_PS10} S.~Papadima, A.I.~Suciu,
{\em Bieri--{N}eumann--{S}trebel--{R}enz invariants and 
homology jumping loci}, Proc.~London Math.~Soc. 
\textbf{100} (2010), no.~3, 795--834.
\MRh{2640291}{2011i:55006}

\bibitem{DS_Ra11} R.~Randell, 
{\em The topology of hyperplane arrangements}, in: 
Topology of algebraic varieties and singularities, 309--318,
Contemp. Math., vol.~538, Amer. Math. Soc., Providence, RI, 2011. 
\MRh{2777827}{2012e:32046}

\end{thebibliography}
\end{document}